\documentclass[14pt, reqno]{amsart}
\usepackage{graphicx,psfrag}
\usepackage[psamsfonts]{amssymb}
\usepackage{pdfsync}
\usepackage{amscd}
\usepackage{amsmath}
\usepackage{mathrsfs}
\usepackage{stmaryrd}
\usepackage[small,nohug,heads=littlevee]{diagrams}
\usepackage{comment}
\usepackage{enumerate}
\usepackage{xfrac}
\usepackage{cases}
\diagramstyle[labelstyle=\scriptstyle]

\theoremstyle{plain}
    \newtheorem{theorem}{Theorem}[section]
    \newtheorem{lemma}[theorem]{Lemma} 
    \newtheorem{corollary}[theorem]{Corollary}
    \newtheorem{proposition}[theorem]{Proposition}
    
    \theoremstyle{definition}
    \newtheorem{definition}[theorem]{Definition}
    
    \newtheorem{example}[theorem]{Example}

    \newtheorem{remark}[theorem]{Remark}

\theoremstyle{remark}

\numberwithin{equation}{section}

    \newcommand{\R}{\mathbb{R}}
    
    \newcommand{\C}{\mathbb{C}} 
    
    \newcommand{\Z}{\mathbb{Z}}

        \newcommand{\g}{\mathfrak{g}} 
   
  \newcommand{\A}{\mathcal{A}}

     \newcommand{\pt}{\text{\it pt}}

\title{Formal Verlinde Module}

\author{Yanli Song}

\date{\today}

\begin{document}
\maketitle
\begin{abstract}
Let $G$ be a compact, simple and simply connected Lie group and $\A$ be an equivariant Dixmier-Douady bundle over $G$. For any fixed level $k$, we can define a $G$-$C^{*}$-algebra $C_{\A^{k+h^{\vee}}}(G)$ as all the continuous sections of the tensor power $\A^{k+h^{\vee}}$ vanishing at infinity.  A deep theorem by Freed-Hopkins-Teleman showed that the twisted $K$-homology $KK_{\bullet}^{G}(C_{\A^{k+h^{\vee}}}(G), \C)$ is isomorphic to the level $k$ Verlinde ring $R_{k}(G)$. By the construction of crossed product, we define a $C^{*}$-algebra $C^{*}(G,C_{\A^{k+h^{\vee}}}(G))$. We show that the $K$-homology $KK_{\bullet}(C^{*}(G,C_{\A^{k+h^{\vee}}}(G)),\C)$ is isomorphic to the formal Verlinde module $R^{-\infty}(G) \otimes_{R(G)} R_{k}(G)$, where $R^{-\infty}(G) = \mathrm{Hom}_{\Z}(R(G),\Z)$ is the completion of the representation ring. 
\end{abstract}

\tableofcontents

\abovedisplayskip=2pt
\belowdisplayskip=2pt

\section{Introduction}

Let $G$ be a compact, simple, and simply connected  Lie group. The level $k$ Verlinde ring $R_{k}(G)$ is defined as the ring of positive energy representations of the loop group $LG$ at level $k$. We can also define $R_{k}(G) = R(G) /I_{k}(G)$, where $I_{k}(G)$ is the level $k$ fusion ideal. 

On the other hand, one can construct a $G$-equivariant Dixmier-Douady bundle $\A \to G$ corresponding to the generator 
\[
[1] \in H^{3}_{G}(G) \cong \Z.
\] 
 All the continuous sections of its tensor power $\A^{k+h^{\vee}}$($h^{\vee}$ is the dual Coxeter number of $G$) vanishing at infinity  form a $G$-$C^{*}$-algebra $C_{\A^{k+h^{\vee}}}(G)$. With this $C^{*}$-algebra, we can define the equivariant twisted $K$-homology 
\[
KK^{G}_{\bullet}(C_{\A^{k+ h^{\vee}}}(G), \C).
\]
Freed-Hopkins-Teleman \cite{Freed11,Freed11-2,Freed13} showed that
\begin{equation}
\label{FHT}
KK^{G}_{\bullet}(C_{\A^{k+ h^{\vee}}}(G), \C) \cong R_{k}(G). 
\end{equation}
In fact, they proved the isomorphism (\ref{FHT}) for more general Lie groups. 

For a special case when $G$ is compact, simple, and simply connected,   Kitchloo-Morava \cite{Kitchloo07} and  Meinrenken \cite{Meinrenken09})  developed different approaches to (\ref{FHT}). This paper is based on Meinrenken's method. We first realize geometrical the Dixmier-Douady bundle $\A \to G$ as a co-simplicial bundle. In more detail,  suppose that $\{0, \dots, l = \mathrm{dim}(G)\}$ label the vertices of the fundamental alcove $\Delta$ of $G$. For any $I \subseteq \{0, \dots, l\}$, let $\Delta_{I}$ be the face of $\Delta$ spanned by vertices in $I$ and $G_{I}$ be the corresponding closed sub-group of $G$.  The simplices of $\A$ are given by
\[
\A_{I} = G \times_{G_{I}} \mathcal{K}(H)  \to G / G_{I}, 
\] 
where  $\mathcal{K}(H)$ denotes the compact operators on Hilbert space $H$. In addition, we can associate each sub-group $G_{I}$ with a central extension $\hat{G}_{I}$, which is determined by the Dixmier-Douady bundle $\A$. 

When we restrict to the simplices, the equivariant twisted $K$-homology 
\[
KK^{G}_{\bullet}(C_{\A^{k+ h^{\vee}}}(G/G_{I}), \C) \cong KK^{G_{I}}_{\bullet}(\mathcal{K}(H^{k}), \C)
\]
is isomorphic to $R(\hat{G}_{I})_{k}$, which is the Grothendieck group of all the $\hat{G}_{I}$-representations on which the central $U(1)$ acts with weight $k$. Hence we obtain a spectral sequence consisting of different representation rings. Moreover,  the spectral sequence collapses at the $E^{2}$-stage by Bott-periodicity. All these provide us a way of computing the equivariant twisted $K$-homology $KK^{G}_{\bullet}(C_{\A^{k+ h^{\vee}}}(G), \C)$ in terms of weights.  

In this paper, we extend the identity (\ref{FHT}) to the infinite case. To be more precise, since $C_{\A^{k+h^{\vee}}}(G)$ is a $G$-$C^{*}$-algebra, we can construct its crossed product $C^{*}(G, C_{\A^{k+h^{\vee}}}(G))$, which is a $C^{*}$-algebra. Let us consider the $K$-homology
\[
KK_{\bullet}(C^{*}(G, C_{\A^{k+h^{\vee}}}(G)), \C).
\]
By the simplicial description of the Dixmier-Douady bundle $\A$, we can restrict to the simplices as well. By some basic properties of $C^{*}$-algebras and induced representation of $C^{*}$-algebras, we can show that
\[
KK_{\bullet}(C^{*}(G, C_{\A^{k+h^{\vee}}_{I}}(G/G_{I})), \C) \cong KK_{\bullet}(C^{*}(G_{I}, \mathcal{K}(H^{k})), \C),  
\]
and the right-hand side of the above equation is isomorphic to 
\[
R^{-\infty}(\hat{G}_{I})_{k} = \mathrm{Hom}_{\Z}(R(\hat{G}_{I})_{k},\Z),
\] 
the completion of  the representation ring  $R(\hat{G}_{I})_{k}$. As in the classical case, by constructing a spectral sequence, we can compute the $K$-homology $KK_{\bullet}(C^{*}(G, C_{\A^{k+h^{\vee}}}(G)), \C)$ in terms of weights. To summarize, 
\begin{theorem}[Main Theorem]
Let $G$ be a compact, simple, and simply connected Lie group. We have the isomorphism:
\begin{equation}
KK_{\bullet}(C^{*}(G, C_{\A^{k+h^{\vee}}}(G)), \C) \cong R^{-\infty}(G) \otimes_{R(G)} R_{k}(G), 
\end{equation}
where $R^{-\infty}(G) = \mathrm{Hom}_{\Z}(R(G), \Z)$, is the completion of the character ring $R(G)$. 
\end{theorem}

Our motivation to consider the $K$-homology for crossed products is a undergoing project of understanding the quantization of quasi-Hamiltonian spaces.  We will explain this idea in the last section. 

$\mathbf{Acknowledgments}$: The author would like to thank Nigel Higson and Eckhard Meinrenken for many benefited discussions and useful advices.


\section{Dixmier-Douady Bundles}
In this section we will  review the  Dixmier-Douady theory \cite{Dixmier63} and  give a simplicial construction of the Dixmier-Douady bundle $\A$ over the Lie group $G$. 

Let $X$ be a CW-complex, and $H$ a separable Hilbert space. We denote by $U(H), \mathcal{K}(H)$ the unitary operators and compact operators on $H$ respectively. There is an isomorphism
\[
\mathrm{Aut}(\mathcal{K}(H)) \cong PU(H), 
\]
where $PU(H)$ is the projective unitary group. 

\begin{definition}[\cite{Dixmier63, Rosenberg89,Williams98}]
A \emph{Dixmier-Douady} bundle over $X$ is a locally trivial bundle $\A \to X$ of $C^{*}$-algebras, with typical fiber $\mathcal{K}(H)$, and  structure group $PU(H)$.  That is, 
\[
\A = \mathcal{P} \times_{PU(H)} \mathcal{K}(H),
\]
where $\mathcal{P} \to X$ is a principle $PU(H)$-bundle over $X$. 
\end{definition}

\begin{example}
Let $V \to X$ be an even-dimensional Euclidean vector bundle and $\mathrm{Cliff}(V) \to X$ the complex Clifford algebra bundle over $X$. For each $x \in X$, $\mathrm{Cliff}(V_{x})$ is a matrix algebra. Hence, $\mathrm{Cliff}(V)$ is an example of Dixmier-Douady bundle with finite rank. 
\end{example}

If $\A_{1}, \A_{2}$ are two Dixmier-Douady bundles over $X$, with fibers isomorphic to $\mathcal{K}(H_{1})$ and $\mathcal{K}(H_{2})$, then their tensor product $\A_{1}\otimes \A_{2}$ is a Dixmier-Douady bundle with fiber modeled on $\mathcal{K}(H_{1} \otimes H_{2})$. In addition, the opposite bundle $\A^{\mathrm{opp}}$ is a Dixmier-Douady bundle with fiber modeled on $\mathcal{K}(\overline{H})$, where $\overline{H}$ denotes the conjugate Hilbert space of $H$.

\begin{definition}
A \emph{Morita trivialization} of a Dixmier-Douady bundle $\A \to X$ is given by a bundle of Hilbert space $\mathcal{E} \to X$ and an isomorphism $\psi : \mathcal{K}(\mathcal{E}) \cong \A$.  
We say that two Dixmier-Douady bundles $\A, \mathcal{B}$ are \emph{Morita equivalent} if there exists a Morita trivialization for $\A \otimes \mathcal{B}^{\mathrm{opp}}$.
\end{definition}

If $M$ is an even-dimensional, smooth manifold, then $\mathrm{Cliff}(TM) \to M$ is a Dixmier-Douady bundle over $M$. 
As observed by Connes \cite{connes85} and Plymen \cite{plymen86}, a Spin$^{c}$-structure for $M$ can be defined as a Morita trivialization:
\[
\mathcal{S} : \mathrm{Cliff}(TM) \sim \C. 
\]
There is a classical result of Dixmier-Douady bundles. 

\begin{theorem}[\cite{Dixmier63}]
The Morita equivalence classes of Dixmier-Douady bundles are isomorphic to $H^{3}(X)$. 
\end{theorem}

\begin{definition}
\label{morita-1}
Suppose that $\A \to X, \mathcal{B} \to Y$ are two Dixmier-Douady bundles modeled on $\mathcal{K}(H_{A})$ and $\mathcal{K}(H_{B})$. A \emph{Morita morphism}
\[
(\phi, \mathcal{E}) : \A \to \mathcal{B}
\]
is a proper map $\phi : X \to Y$ together with a bundle of Banach space $\mathcal{E} \to X$, which carries a bi-module structure:
\[
\phi^{*}(\mathcal{B}) \circlearrowleft \mathcal{E} \circlearrowright \A,
\]
locally modeled on $\mathcal{K}(H_{B})  \circlearrowleft \mathcal{K}(H_{B}, H_{A})  \circlearrowright \mathcal{K}(H_{A})$.
\end{definition}

\begin{remark}
\label{2-iso}
Suppose that $\mathcal{E}_{1}, \mathcal{E}_{2}$ are two Morita morphisms between $\A$ and $\mathcal{B}$. They are related by a line bundle $L \to X$:
\[
L := \mathrm{Hom}_{\phi^{*}(\mathcal{B})-\A}(\mathcal{E}_{1}, \mathcal{E}_{2}).
\]
\end{remark}

Let now $G$ be a compact,  simple, and simply connected Lie group. All of above extend to the equivariant setting by straightforward modification. In particular, the Morita equivalence classes of $G$-equivariant Dixmier-Douady bundles are classified by the third equivariant cohomology $H^{3}_{G}(X)$.

\begin{lemma}
\label{central extension}
There is an 1-1 correspondence between $G$-equivariant Dixmier-Douady bundles over a point and central extensions of $G$.
\begin{proof}
For a $G$-equivariant Dixmier-Douady bundle $\A \to \pt$, we can obtain a central extension as a pull-back of $U(H) \to \mathrm{Aut}(\mathcal{K}(H))$ by the homomorphism $G \to \mathrm{Aut}(\mathcal{K}(H))$. 

On the other hand, let $\hat{G}$ be an arbitrary central extension of $G$. We can define a Hilbert space $H$ to be all the $L^{2}$-sections of the associated line bundle $\hat{G} \times_{U(1)} \C$. Then
\[
\A = \mathcal{K}(H) \to \pt
\] 
defines a  $G$-equivariant Dixmier-Douady bundle over a point. 
\end{proof}
\end{lemma}

By the assumptions on the Lie group $G$, we have that $H^{3}_{G}(G) \cong \Z$. Let us focus on the construction of the specific Diximer-Douady bundle $\A$ which corresponds to the generator $[1] \in H^{3}_{G}(G)$.

Let $T \subset G$ be a maximal torus of $G$ with Lie algebra $\mathfrak{t}$. We choose a positive Weyl chamber $\mathfrak{t}_{+}$ and a closed fundamental  Weyl alcove $\Delta$. We label the vertices of $\Delta$ by $0, \dots, l = \mathrm{rank}(G)$ such that $0$ corresponds to the origin. For every non-empty subset $I \subseteq \{0,\dots l\}$, it uniquely determines a closed wall $\Delta_{I}$ in $\Delta$ which is the simplex spanned by vertices in $I$. For any point $x$  in the interior of $\Delta_{I}$, we denote by $G_{I}$ the isotropy group of $\mathrm{exp}(x) \in G$ associated to the conjugate $G$-action. The group $G_{I}$ does not depend on the choices of point $x$. In particular, we have that
\[
G_{\{0\}} = G, \ G_{\{0, \dots, l\}} = T. 
\]

\begin{theorem}
There is an isomorphism:
\begin{equation}
\label{decomposition of G}
G \cong \coprod_{I}G/G_{I} \times \Delta_{I} / \sim,
\end{equation}
where the gluing identification are given by
\begin{equation}
\label{glue}
(g, \iota_{J}^{I}(\xi)) \sim (\phi_{I}^{J}(g), \xi), J \subset I. 
\end{equation}
Here $\iota_{J}^{I} : \Delta_{J} \hookrightarrow \Delta_{I}$ is the inclusion and $\phi_{I}^{J} : G/G_{I} \to G/G_{J}$ is the natural projection.
\begin{proof}
\cite[section 3]{Meinrenken09}.
\end{proof}
\end{theorem}

With the description of $G$ in (\ref{decomposition of G}), we can associate the Dixmier-Douady bundle $\A$ with a $G$-equivariant collection of central extensions \cite[Lemma 3.3]{Meinrenken09}: 
\begin{equation}
\label{central extension} 
1 \to U(1) \to \hat{G}_{I} \to G_{I} \to 1, 
\end{equation}
together with lifts $\hat{\iota}_{I}^{J} : \hat{G}_{I} \hookrightarrow \hat{G}_{J}$ of the inclusions $\iota_{I}^{J} : G_{I} \hookrightarrow G_{J}$ for $J \subset I$, satisfying the coherence condition:
\[
\hat{\iota}_{I}^{K}  = \hat{\iota}_{J}^{K} \circ \hat{\iota}_{I}^{J}  \ \mathrm{for} \ K \subset J \subset I.
\]

\begin{example}
If $G= SU(2)$, then $\Delta = [0,1]$. There are only two vertices 0 and 1. The central extensions $\hat{G}_{0} \cong \hat{G}_{1} \cong G \times U(1)$ and $\hat{G}_{\{0,1\}} \cong T \times U(1)$. The inclusion  of $\hat{G}_{\{0,1\}}$ into $\hat{G}_{0}$ is given by $(t,z) \mapsto (t,z)$ while the inclusion  of $\hat{G}_{\{0,1\}}$ into $\hat{G}_{1}$ is given by $(t,z) \mapsto (t, t^{\rho}\cdot z)$. 
\end{example}

In addition, by the theory of affine Lie algebras, we can construct a Hilbert space $H$, equipped with unitary representation of the central extensions $\hat{G}_{I}$ such that the central $U(1)$ acts with weight -1 and the action of $\hat{G}_{J}$ restricts to the action of $\hat{G}_{I}$ for any $J \subset I$ \cite[Lemma 3.5]{Meinrenken09}. Let $\A_{I} = G \times_{G_{I}} \mathcal{K}(H)$. For any $J \subset I$, there is a homomorphism $\A_{I} \to \A_{J}$ coving the map $\phi_{I}^{J} : G/G_{I} \to G/G_{J}$. We can glue all the $\A_{I}$ together so as to recover the $G$-equivariant Dixmier-Douady bundle over $G$:
\begin{equation}
\label{G}
\A = \coprod_{I}(\A_{I} \times \Delta_{I}) / \sim,
\end{equation}
where the gluing identifications come from (\ref{glue}). Here (\ref{G}) indeed gives a simplicial description of the Dixmier-Douady bundle $\A$.
 
There is an alternative construction \cite{Atiyah04} of the Dixmier-Douady bundle $\A \to G$:
\[
\A = P_{e}G \times_{L_{e}G} \mathcal{K}(H),
\]
where $P_{e}G$ is the spaces of based paths in $G$, $L_{e}G$ is the based loop group, and $H$ is a representation of the standard central extension $\widehat{LG}$ such that the central circle acts with weight -1.


\section{The Morita Categories}
In this section we want to give some theoretical background for the discussion of $C^{*}$-algebras, crossed products and induced representations of $C^{*}$-algebras.  The systematic approach to all this in terms of Hilbert modules are based on the fundamental idea of Rieffel \cite{Rieffel74}. We use \cite{Williams07} as our primary reference.
 
\begin{definition}
If $B$ is a $C^{*}$-algebra, then a (right) \emph{pre-Hilbert $B$-module} is a complex Banach space $E$ equipped with a (right) $B$-action and a $B$-valued inner product $\langle \  , \  \rangle_{B} : E \times E \to B$, which is linear in the second variable and anti-linear in the first variable and satisfies 
\begin{itemize}
\item
$(\langle \xi, \eta \rangle_{B})^{*} = \langle \eta, \xi \rangle_{B}$;
\item
$\langle \eta, \xi \rangle_{B} b= \langle \eta, \xi  \cdot b \rangle_{B}$;
\item
$\| \xi \|^{2} = \| \langle \xi, \xi \rangle_{B}\|$
\end{itemize}
for all $\xi, \eta \in E, b \in B$. The norm completion of $E$ is a \emph{Hilbert $B$-module}.
\end{definition}

The Hilbert $\C$-modules are precisely the Hilbert spaces. 

\begin{definition}

Let $E$ be a Hilbert $B$-module.  We say that a linear map $T : E \to E$ is \emph{adjointable} if there exists a map $T^{*} : E \to E$ such that 
\[
\langle T\xi, \eta \rangle_{B} = \langle \xi, T^{*} \eta \rangle_{B}.
\] 
We denote $\mathcal{L}_{B}(E) =  \{T : E \to E : T \ \mathrm{is \ adjointable} \}$.
\end{definition}

Let $A, B$ be two $C^{*}$-algebras. 

\begin{definition}
A \emph{Hilbert A-B bi-module} is a pair $(E, \phi)$ in which $E$ is a Hilbert $B$-module and $\phi : A \to \mathcal{L}_{B}(E)$ is a $*$-representation of $A$ on $E$. Two Hilbert $A$-$B$ bi-modules $(E_{i}, \phi_{i}), i = 1, 2$ are called \emph{unitarily equivalent} if there is an isomorphism $U : E_{1} \to E_{2}$ preserving the $B$-valued inner product and commuting with $\phi$.  In addition, we say that $(E, \phi)$ is \emph{non-degenerate} if $\phi(A)E = E$. \emph{Morita morphisms} between $A$ and $B$ are given by non-degenerate Hilbert $A$-$B$ bi-modules. 
\end{definition}


\begin{definition}
The \emph{Morita category} $\mathfrak{M}$ is the category whose objects are $C^{*}$-algebras and the morphisms are given by  unitary equivalence classes of non-degenerate Hilbert bi-modules. 
\end{definition}

If $G$ is a compact Lie group, we can define the $G$-equivariant Morita category in a similar way. 

\begin{definition}
A  $G$-$C^{*}$-algebra $A$ is a $C^{*}$-algebra equipped with a $*$-homomorphism 
\[
\alpha : G \to \mathrm{Aut}(A): s \mapsto \alpha_{s},
\]
such that $s \mapsto \alpha_{s}$ is strongly continuous. 
\end{definition}

Let $A,B$ be two $G$-$C^{*}$-algebras. 

\begin{definition}
A \emph{G-equivariant Hilbert A-B bi-module} is given by $(E, \phi, u)$ where $(E, \phi)$ is a non-equivariant Hilbert $A$-$B$ bimodule, and $u$ is a strongly continuous homomorphism $u : G \to \mathrm{Aut}(E)$ such that
\begin{itemize}
\item
$\langle u_{g}(\xi), u_{g}(\eta) \rangle_{B} = g \cdot \langle \xi, \eta \rangle_{B}$;
\item
$u_{g}(\xi \cdot b) = u_{g}(\xi)\cdot (g \cdot b)$;
\item
$u_{g}(\phi(a) \xi) = \phi(g \cdot a) u_{g}(\xi)$.
\end{itemize}
\end{definition}

\begin{definition}
The \emph{$G$-equivariant Morita category} $\mathfrak{M}(G)$ is the category whose objects are $G$-$C^{*}$-algebras and the morphisms are given by unitary equivalence classes of non-degenerate $G$-equivariant Hilbert bi-modules. 
\end{definition}

There is a descent operation which transfers the equivariant Morita category  to the non-equivariant Morita category. It involves the notion of crossed products of $C^{*}$-algebras. We will give a rapid review on the basic construction.

If $A$ is $G$-$C^{*}$-algebra, then the space $C(G, A)$ of continuous $A$-valued functions on $G$ becomes a $*$-algebra with respect to the convolution and involution defined by
\[
f * g(s) = \int_{G} f(t) \alpha_{t}(g(t^{-1}s))dt,
\]
and
\[
f^{*}(s)  = \alpha_{s}(f(s^{-1}))^{*}. 
\]
There is an injection $C(G, A) \hookrightarrow \mathcal{B}(L^{2}(G, A))$:
\[
f \mapsto T_{f}, T_{f}(g) = f *g. 
\]
\begin{definition}
We define the \emph{crossed product} $C^{*}(G, A)$ as the completion of $C(G, A)$ in $\mathcal{B}(L^{2}(G, A))$ with respect to the norm $\|f\| = \|T_{f}\|$.
\end{definition}

\begin{remark}
The $C^{*}$-algebra $C^{*}(G, A)$ is indeed the so-called reduced crossed product. There is an alternative definition of crossed product, namely the full crossed product, which is the completion of $C(G, A)$ with respect to different norms. However, since $G$ is  compact, thus \emph{amenable}. The two crossed products coincide. Thus, we will not distinguish them in this paper.  
\end{remark}
If $A = \C$, then we get the group $C^{*}$-algebra $C^{*}(G)$. There is a Morita equivalence: 
\[
C^{*}(G) \cong \bigoplus_{\sigma \in \mathrm{Irrep}(G)}M_{n}(\C),
\]
where $M_{n}(\C)$ is the matrix algebra with $n = \mathrm{dim}(\sigma)$.

We can extend the construction of crossed products of $C^{*}$-algebras to the morphisms. Let $(E, \phi)$ be a Hilbert $A$-$B$ bi-module, equipped with a strongly continuous homomorphism $u : G \to \mathrm{Aut}(E)$ such that
\[
\langle u_{s}(\xi) , u_{s}(\eta) \rangle_{B} = s \cdot (\langle \xi, \eta \rangle_{B}),
\]
and
\[
 u_{s}(\xi \cdot b) = u_{s}(\xi) \cdot (s \cdot b), \  u_{s}(\phi(a) \xi ) = \phi(s \cdot a) \cdot u_{s}(\xi),
\]
for $\xi, \eta \in E, a \in A, b \in B, s \in G$. The crossed product $(\hat{E}, \hat{\phi})$ is a Hilbert $C^{*}(G,A)$-$C^{*}(G,B)$ bi-module, defined as the completion of $C(G, E)$ with respect to the $C^{*}(G,B)$-valued inner product:
\[
\langle \xi, \eta \rangle (t) = \int_{G} \langle \xi(s), u_{s} (\eta(s^{-1}t))\rangle ds;
\]
and with left action of $C(G, A)$ on $C(G, E)$ given by
\[
\hat{\phi}(f)(\xi)(t) = \int_{G} \phi(f(s)) u_{s}(\xi(s^{-1}t))ds.
\]
The Hilbert bi-module $(\hat{E},\hat{\phi})$ gives a Morita morphism between crossed products $C^{*}(G,A)$ and $C^{*}(G,B)$ 

If $H \subseteq G$ is a closed subgroup of $G$, then we can define  an induction  functor:
\[ 
\mathcal{I}_{H}^{G}: \mathfrak{M}(H) \to \mathfrak{M}(G).
\]
 
\begin{definition}
Let $A$ be a $H$-$C^{*}$-algebra. The induced $G$-$C^{*}$-algebra $\mathcal{I}_{H}^{G}(A)$ is defined as
\[
 \{ f \in C(G, A): f(sh) = u_{h^{-1}}( f(s)), \ \mathrm{and} \ sH \mapsto \|f(s)\| \in C(G/H) \},
\] 
equipped with the point-wise operations and the supremum norm. The induced $G$-action on $\mathcal{I}_{H}^{G}(A)$ is given by
\[
u_{s}(f)(t) = f(s^{-1}t), \ \mathrm{for} \ s,t \in G.
\]
\end{definition}

\begin{example}
if $A = C_{0}(X)$ for some topological $H$-space $X$, then 
\[
\mathcal{I}_{H}^{G}(A) = C_{0}(G \times_{H}X).
\]
\end{example}
Suppose that $(E, \phi)$ is a $H$-equivariant Morita morphism between two $H$-$C^{*}$-algebras $A$ and $B$. An obvious extension of the above construction yields an induced Morita morphism $\mathcal{I}_{H}^{G}(E, \phi)$, which is a non-degenerate Hilbert $\mathcal{I}_{H}^{G}(A)$-$\mathcal{I}_{H}^{G}(B)$ bi-module. The Green's imprimitivity theorem  \cite{Green78} reads as follows:
\begin{theorem}
Let $G$ be a compact Lie group and $H \subseteq G$ a closed subgroup. Suppose that $A \in \mathfrak{M}(H)$ is a $H$-$C^{*}$-algebra. Its induced $C^{*}$-algebra $\mathcal{I}_{H}^{G}(A) \in \mathfrak{M}(G)$. We have that
\[
C^{*}(H, A) \ \mathrm{and} \ C^{*}(G, \mathcal{I}_{H}^{G}(A)) 
\]
are Morita equivalent in $\mathfrak{M}$.
\end{theorem}

\section{Equivariant Twisted K-Homology}
\label{twisted K-homology}
Suggested by Rosenberg's definition of twisted $K$-theory \cite{Rosenberg89}, we define twisted $K$-homology  as the $K$-homology of the $C^{*}$-algebra of sections of the Dixmier-Douady bundle $\A$. We begin by reviewing  Kasparov's approach to the analytic $K$-homology. We refer to Kasparov's papers \cite{Kasparov80, Kasparov88} and Higson-Roe's book \cite{Higson00}. 

\begin{definition}
Let $A$ be a separable C$^{*}$-algebra. An (even) \emph{Fredholm module} is a triple consisting of
\begin{itemize}
\item
 a separable $\Z_{2}$-graded Hilbert space $H$;
\item
 an even $*$-homomorphism:  $\rho : A \to \mathcal{B}(H)$;
\item
 an odd bounded operator $F$ on $H$ such that for all $a \in A$
\[
(F - F^{*}) \rho(a) \in \mathcal{K}(H), \ (F^{2} - 1) \rho(a) \in \mathcal{K}(H), \  [\rho(a), F]  \in \mathcal{K}(H).
\]
Here $\mathcal{K}(H)$  and $\mathcal{B}(H)$ denote the compact operators and the bounded operators on $H$ respectively. 
\end{itemize}

\end{definition}

An \emph{operator homotopy} between two Fredholm modules $(H, F_{0}, \rho)$ and $(H, F_{1}, \rho)$ is a norm continuous path of Fredholm modules
\[
t \mapsto (H, F_{t}, \rho), \ t \in [0,1].
\]
Two Fredholm modules are \emph{equivalent} if they are related by \emph{unitary transformations} or \emph{operator homotopies}.

\begin{definition}
The set of equivalence classes of Fredholm modules $(H, F, \rho)$ forms an abelian group $KK_{0}(A,\C )$, where the addition is given by direct sum and the \emph{zero module} has zero Hilbert space, zero representation, and zero operator.   
\end{definition}

More general, one defines $KK_{n}(A, \C) = KK_{0}(\mathrm{Cliff}(\R^{n})\otimes A, \C)$ and it has the mod 2 periodicity property: $KK_{n}(A, \C) = KK_{n+2}(A, \C)$. Thus,  $KK_{\bullet}(A, \C)$ is a $\Z_{2}$-graded homology theory on $C^{*}$-algebras. In a special case when $A = \C$, we have that
\[
KK_{0}(\C, \C) \cong \Z, \ KK_{1}(\C, \C) = 0. 
\]
The $K$-homology groups $KK_{\bullet}(A, \C)$ are functorial with respect to \emph{Morita morphisms} between $C^{*}$-algebras. 

Let $G$ be a compact Lie group. For any $G$-$C^{*}$-algebra $A$, we can define the $G$-equivariant $K$-homology $KK^{G}_{\bullet}(A, \C)$, using equivariant Fredholm modules with straightforward modification. 
\begin{example}
If $A = \C$ with trivial $G$-action, we have that
\[
KK^{G}_{0}(\C, \C) \cong R(G), \ KK^{G}_{1}(\C, \C) = 0. 
\]
For the $K$-homology of group $C^{*}$-algebras, we have that
\[
KK_{0}(C^{*}(G),\C) \cong \hat{R}(G), \ KK_{1}(C^{*}(G),\C) =0,
\]
where $\hat{R}(G) = \mathrm{Hom}_{\Z}(R(G), \Z)$ is a completion of the representation ring $R(G)$.
\end{example}


\begin{definition} 
Given a $G$-equivariant Dixmier-Douady bundle $\A$ over a $G$-space $X$, we define a $G$-$C^{*}$-algebra 
\[
C_{\A}(X) = \Gamma_{0}(X, \A) 
\]
as all the sections of $\A$ vanishing at infinity with supremum norm. The \emph{equivariant twisted K-homology} $K^{G}_{\A}(X)$ is defined by $KK^{G}_{\bullet}(C_{\A}(X), \C)$. 
\end{definition}

We list some basic properties of  equivariant twisted $K$-homology, which can be deduced from the equivariant $KK$-theory of general $C^{*}$-algebras \cite[Section 2]{Kasparov88}. 

\begin{lemma}
\label{mortia map}
Suppose that $\A \to X, \mathcal{B} \to Y$ are two $G$-equivariant Dixmier-Douady bundles. If $(\phi, \mathcal{E}) : \A \to \mathcal{B}$ is a Morita morphism, then $(\phi, \mathcal{E})$ induces a map
\[
\phi_{\mathcal{E}} : KK^{G}_{\bullet}(C_{\A}(X), \C) \to KK^{G}_{\bullet}(C_{\mathcal{B}}(Y), \C).
\]
\begin{proof}
The Morita morphism $(\phi, \mathcal{E})$ from $\A$ to $\mathcal{B}$ is indeed a Hilbert $C_{\mathcal{B}}(X)$-$C_{\A}(Y)$ bi-module, that is a Morita morphsim between the two $C^{*}$-algebras $C_{\A}(X)$ and $C_{\mathcal{B}}(Y)$. 
\end{proof}
\end{lemma}

\begin{remark}
As in Remark \ref{2-iso}, if $(\phi, \mathcal{E}_{1})$ and $(\phi, \mathcal{E}_{2})$ are two Morita morphisms between $\A$ and $\mathcal{B}$, then the two maps:
\[
\phi_{\mathcal{E}_{1}} \ \mathrm{and} \ \phi_{\mathcal{E}_{2}}
\]
may not be identical unless the line bundle 
\[
L := \mathrm{Hom}_{\phi^{*}(\mathcal{B})-\A}(\mathcal{E}_{1}, \mathcal{E}_{2}).
\]
can be trivialized. 

\end{remark}

\begin{lemma}
\label{sequence}
Suppose that $X$ is a $G$-$CW$-complex and $Y \subset X$ is a closed $G$-invariant subset. There is a long exact sequence:
\begin{equation}
\begin{aligned} 
&\dots \to KK^{G}_{n}(C_{\A}(Y), \C) \to KK^{G}_{n}(C_{\A}(X), \C) \to\\
&   KK^{G}_{n}(C_{\A}(X \setminus Y), \C) \to   KK^{G}_{n-1}(C_{\A}(Y), \C) \to \dots
\end{aligned}
\end{equation}
Hence, any filtration of $X$ induces a spectral sequence for $KK^{G}_{\bullet}(C_{\A}(X), \C)$. 

\end{lemma}

\begin{lemma}
\label{point}
If $X = \pt$, then $\A \cong \mathcal{K}(H)$. Let $\hat{G}$ be the central extension induced by 
\begin{diagram}
 \hat{G} & \rTo  & U(H)  \\
 \dTo&  & \dTo \\
 G & \rTo  & \mathrm{Aut}(\mathcal{K}(H)) \cong PU(H) \\
\end{diagram}
We have that
\[ 
KK^{G}_{0} (C_{\A}(\pt), \C) \cong R(\hat{G})_{-1}, 
\]
where $R(\hat{G})_{-1}$ is the Grothendieck group of $\hat{G}$-representations on which the central $U(1)$ acts with weight $-1$. 
\begin{proof}
\cite[Appendix B]{Meinrenken09}
\end{proof}
\end{lemma}

\begin{lemma}
\label{induction}
Suppose that $H \subset G$ is a closed sub-group with maximal rank, and  $\mathcal{B} \to \pt$ is a $H$-equivariant Dixmier-Douady bundle. We can form
\[
 \A = G\times_{H} \mathcal{B} \to G/H, 
 \]
and
\[
C_{\g} = G \times_{H} \mathrm{Cliff}(\g/\mathfrak{h}) \to G/H,
\]
 which are two $G$-equivariant Dixmier-Douady bundles over $G/H$. There is an isomorphism:
 \[ 
 \mathcal{I}_{H}^{G} :  KK^{H}_{0}(C_{\mathcal{B}}(\pt), \C) \cong KK^{G}_{0}(C_{\A \otimes C_{\g}}(G/H),\C).
 \]
\begin{proof}
Remember that all the $H$-equivariant Dixmier-Douady bundles over a point are Morita equivalent to finite rank ones. For an arbitrary element in $KK^{H}_{0}(C_{\mathcal{B}}(\pt), \C)$, we assume that it can be realized by a Fredholm module $(V, \rho, 0)$, where $V$ is a finite dimensional Hilbert $H$-space. We define
\[
 \mathcal{E} = G \times_{H} V, \  \  \mathcal{S} = \Lambda^{*}T^{*}(G/H),
\]
and a Hilbert space
\[
\mathcal{H} = L^{2}(G/H, \mathcal{E} \otimes \mathcal{S}).
\]
In addition, let $d$ be the de Rham differential operator, $d^{*}$ its dual, and $\mathcal{D} = d + d^{*}$. Then $\mathcal{D}$ is essentially self-adjoint and
\[
F = \frac{\mathcal{D}}{\sqrt{\mathcal{D}^{2}+1}}
\]  
is a bounded operator on $\mathcal{H}$. The pair $(\mathcal{H},  F)$ constitutes a Fredholm module over $C_{\mathcal{A} \otimes C_{\g}}(G/H)$ and it determines an element in 
\[
KK^{G}_{0}(C_{\A \otimes C_{\g}}(G/H), \C).
\] 
This construction realizes the isomorphism $\mathcal{I}_{H}^{G}$. 
\end{proof}

\begin{remark}
In the above construction, we do not require that the homogeneous space $G/H$ has a complex structure or even a Spin$^{c}$-structure. In fact, the Dixmier-Douady bundle $\mathrm{Cliff}(G/H) \to G/H$ yields a so-called twisted Spin$^{c}$-structure on $G/H$. In \cite{Meinrenken-book, Sjamaar13}, they give a more explicit discussion. 
\end{remark}

\end{lemma}

For the $K$-homology of crossed products, Lemma \ref{mortia map}-\ref{induction} can be generalized as follow:
\begin{enumerate}
\item
Morita morphism $(\phi, \mathcal{E}) : \A \to \mathcal{B}$ induces a map 
\[
\hat{\phi}_{\mathcal{E}} : KK_{\bullet}(C^{*}(G, C_{\A}(X)), \C) \to KK_{\bullet}(C^{*}(G, C_{\mathcal{B}}(Y)), \C).
\]
\item
If $Y \subset X$ is a closed $G$-invariant subset, then there is a long exact sequence:
\begin{equation}
\begin{aligned} 
&\dots \to KK_{n}(C^{*}(G, C_{\A}(Y)), \C) \to KK_{n}(C^{*}(G, C_{\A}(X)), \C) \\
& \to KK_{n}(C^{*}(G, C_{\A}(X \setminus Y)), \C) \to  KK_{n-1}(C^{*}(G, C_{\A}(Y)), \C) \to \dots
\end{aligned}
\end{equation}

\item
If $\A \to \pt$ is a $G$-equivariant Dixmier-Douady bundle, then we have that
\begin{equation}
\label{a-2}
KK_{0} (C^{*}(G, C_{\A}(\pt)), \C) \cong R^{-\infty}(\hat{G})_{-1}, 
\end{equation}
where $R^{-\infty}(\hat{G})_{-1}$ is the completion of $R(\hat{G})_{-1}$:
\[ 
R^{-\infty}(\hat{G})_{-1} := \mathrm{Hom}_{\Z}(R(\hat{G})_{-1}, \Z).
\]

\item
If $\mathcal{B} \to \pt$ is a $H$-equivariant Dixmier-Douady bundle, then we can define
\[
\A = G\times_{H} \mathcal{B}, \ \mathrm{and} \  C_{\g} = G \times_{H} \mathrm{Cliff}(\g/\mathfrak{h}).
\]
The isomorphism in Lemma \ref{induction} extends to:
\begin{equation}
\label{a-3}
 \hat{\mathcal{I}}_{H}^{G} : KK_{0}(C^{*}(H,  C_{\mathcal{B}}(\pt)), \C) \cong KK_{0}(C^{*}(G, C_{\A \otimes C_{\g}}(G/H)), \C) .
\end{equation}

\end{enumerate}


\section{The Induction Maps} 
\label{terms of weights}
As we mentioned in Section 2, for any $J \subset I$, there is a homomorphism 
\[
\A_{I} =G \times_{G_{I}} \mathcal{K}(H)   \to A_{J} = G \times_{G_{J}} \mathcal{K}(H)
\] covering the projection map $\phi_{I}^{J} : G/G_{I} \to G/G_{J}$. It induces maps between $K$-homology groups:
\[
\phi_{I,*}^{J} : KK^{G}_{0}(C_{\A_{I}}(G/G_{I}), \C) \to KK^{G}_{0}(C_{\A_{J}}(G/G_{J}), \C),  
\]
and 
\[
\hat{\phi}_{I,*}^{J} : KK_{0}(C^{*}(G, C_{\A_{I}}(G/G_{I})), \C) \to KK_{0}(C^{*}(G, C_{\A_{J}}(G/G_{J})), \C).  
\]
The main goal of this section is to compute the induction maps in terms of weights. 

\begin{theorem}
For any $I \subseteq \{0, \dots, l\}$ and positive integer $k$, 
\[
KK_{0}^{G}(C_{\A^{k+h^{\vee}}_{I}}(G/G_{I}),\C) \cong  R(\hat{G}_{I})_{k},
\]
and
\[
KK_{0}(C^{*}(G, C_{\A^{k+h^{\vee}_{I}}_{I}}(G/G_{I})), \C) \cong R^{-\infty}(\hat{G}_{I})_{k},
\]
where $h^{\vee}$ is the dual Coxeter number, $R(\hat{G}_{I})_{k}$ is the Grothendieck group of $\hat{G_{I}}$-representations where the central $U(1)$ acts with weight $k$,  and $R^{-\infty}(\hat{G}_{I})_{k}$ is the completion of $R(\hat{G}_{I})_{k}$. 
\begin{proof}
For any $I$, Meinrenken \cite[Theorem 4.7]{Meinrenken09} showed that there is a $G_{I}$-equivariant Morita isomorphism
\[
\mathrm{Cliff}(\g/\g_{I}) \cong \mathcal{K}(H^{h^{\vee}}),
\]
where $H$ is the Hilbert space constructed in Section 2. Hence, by Lemma \ref{point} and \ref{induction}, we have that
\begin{equation}
\begin{aligned}
&KK_{0}^{G}(C_{\A^{k+h^{\vee}}}(G/G_{I}),\C)  = KK_{0}^{G}(C_{\A^{k}\otimes C_{\g}}(G/G_{I}),\C) \\
&= KK_{0}^{G_{I}}(C_{\A^{k}}(\pt),\C) = KK_{0}^{G_{I}}(\mathcal{K}(H^{k}),\C)= R(\hat{G}_{I})_{k}.
\end{aligned}
\end{equation}
The proof to the second equation is similar. 
\end{proof}
\end{theorem}

For any $J \subset I$, the homogeneous space $G_{J}/G_{I} = \hat{G}_{J}/\hat{G}_{I}$ has a complex structure on it. Let 
\[
\mathrm{Ind}_{I}^{J} : R(\hat{G}_{I})_{k} \to R(\hat{G}_{J})_{k},
\]
be the holomorphic induction map. The identifications 
\[
KK_{0}^{G}(C_{\A^{k+h^{\vee}}}(G/G_{I}),\C) \cong  R(\hat{G}_{I})_{k}
\] 
in the above theorem intertwine the push-forward maps $\phi_{I,*}^{J}$ with $\mathrm{Ind}_{I}^{J}$. To summarize, 

\begin{theorem}
\label{commutative diagram}
For any $J \subset I$, we have the following commutative diagram:
\begin{diagram}
KK_{0}^{G}(C_{\A^{k+h^{\vee}}}(G/G_{I}),\C) & \rTo^{\cong} &   R(\hat{G}_{I})_{k}  \\
 \dTo^{\phi_{I,*}^{J}} &   & \dTo^{\mathrm{Ind}_{I}^{J}}  \\
KK_{0}^{G}(C_{\A^{k+h^{\vee}}}(G/G_{J}),\C) & \rTo^{\cong}  & R(\hat{G}_{J})_{k}. \\
\end{diagram}
 \begin{proof}
\cite[Section 4]{Meinrenken09}
\end{proof}
\end{theorem}

Correspondingly, in the crossed product case, we have that

\begin{theorem}
\label{commutative diagram-2}
For any $J \subset I$, the diagram commutes
\begin{diagram}
KK_{0}(C^{*}(G, C_{\A^{k+h^{\vee}}_{I}}(G/G_{I})), \C) & \rTo^{\cong} &   R^{-\infty}(\hat{G}_{I})_{k}  \\
 \dTo^{\hat{\phi}_{I,*}^{J}} &   & \dTo^{\mathrm{Ind}_{I}^{J}}  \\
KK_{0}(C^{*}(G, C_{\A^{k+h^{\vee}}_{J}}(G/G_{J})), \C)) & \rTo^{\cong}  & R^{-\infty}(\hat{G}_{J})_{k} \\
\end{diagram}
\end{theorem}

 We now fix a maximal torus $T$ of $G$. Let $W$ be the Weyl group and $\Lambda^{*}$ the weight lattice. Let $\Z[\Lambda^{*}]^{W}$ be the invariant part of $\Z[\Lambda^{*}]$ for the  Weyl group action. We can identify
\[
R(T) \cong \Z[\Lambda^{*}], \ \ R(G) \cong \Z[\Lambda^{*}]^{W}.
\]
Let $\mathfrak{G}_{0} = \{ \alpha_{1}, \dots, \alpha_{l}\}$ be a set of simple roots of $\g$,  $\alpha_{i}^{\vee}$ the corresponding co-roots. We denote by $\alpha_{0} = -\alpha_{\mathrm{max}}$ minus the highest root. The positive Weyl chamber is determined by
\[
\mathfrak{t}_{+}^{*} = \{ \nu \in \mathfrak{t}^{*} \big| \langle \nu, \alpha_{i}^{\vee} \rangle \geq 0, i = 1, \dots, l \}.
\]
We denote by $W_{I}$ the Weyl group of $G_{I}$. In particular, for any $I = \{ i\}$, $W_{i}$ has only two elements $\{e, w_{i}\}$. We denote by $W_{\mathrm{aff}}$ the affine Weyl group generated by $\{w_{i} \}_{i=0}^{l}$.

For the collection of central extensions $\{\hat{G}_{I} \}$ defined in (\ref{central extension}), each  $\hat{G}_{I}$ has maximal torus  $\hat{T} = T \times U(1)$ with weight lattice
\[ 
\hat{\Lambda}^{*} \cong \Lambda^{*} \times \Z \subset \hat{\mathfrak{t}}^{*} = \mathfrak{t}^{*} \times \R.
\]
The simple roots for  $\hat{G}_{I}$ are given by $\mathfrak{G}_{I} = \{ (\alpha_{i}, 0) \}_{i \notin I}$, and its co-roots are consisted of 
\[
(\alpha_{i}^{\vee}, \delta_{i,0}) \in \hat{\mathfrak{t}}= \mathfrak{t} \times \R,  \ i \notin I.
\]
Hence, the positive Weyl chamber of $\hat{G}_{I}$ is given by
\[
\hat{\mathfrak{t}}_{+,I}^{*} = \{ (\nu, s) \in \mathfrak{t}^{*} \times \R \big| \langle \nu, \alpha_{i}^{\vee} \rangle + s \delta_{i,0} \geq 0, i \notin I \}.
\] 
All the irreducible $\hat{G}_{I}$-representations for which the central circle acts with weight $k$ are labeled by $\hat{\mathfrak{t}}_{+,I}^{*} \cap (\Lambda^{*} \times \{ k \})$. 
The  $W_{I}$-action on $\hat{\Lambda}^{*}$ preserves the level and the level-$k$ action is given by
\[ 
w \cdot_{k} \nu  = w(\nu - k \nu_{I}) + k \nu_{I}, \nu \in \Lambda^{*}, 
\]
where $\nu_{I} = \frac{1}{h^{\vee}}(\rho - \rho_{I})$, and $\rho_{I}$ is the half sum of all positive roots of $\g_{I}$.

\begin{definition}
We denote by $\Z[\Lambda^{*}]^{W_{I}}$ the anti-invariant part for the shifted $W_{I}$-action. That is, $v \in \Z[\Lambda^{*}]^{W_{I}}$ if and only if for all $w \in W_{I}$
\[
w \star_{k} v = w \cdot_{k + h^{\vee}}v= (-1)^{\mathrm{length}(w)}\cdot v.
\] 
\end{definition}

\begin{definition}
There is a natural skew-symmetrization map $Sk_{I}: \Z[\Lambda^{*}] \to \Z[\Lambda^{*}]^{W_{I}}$ defined by
\[
Sk_{I} : \nu \mapsto \sum_{w \in W_{I}} (-1)^{\mathrm{length}(w)} w \star_{k} \nu \in \Z[\Lambda^{*}]^{W_{I}}. 
\]
\end{definition}

\begin{theorem}\cite[Lemma 5.1]{Meinrenken09}
\label{identify theorem}
There is an isomorphism 
\[
R(\hat{G}_{I})_{k} \cong \Z[\Lambda^{*}]^{W_{I}}, 
\]
such that it intertwines the holomorphic induction map 
\[ 
\mathrm{Ind}_{I}^{J} : R(\hat{G}_{I})_{k} \to R(\hat{G}_{J})_{k}
\] 
with
\[
Sk_{I}^{J} = \frac{1}{|W_{I}|} Sk^{J} :  \Z[\Lambda^{*}]^{W_{I}} \to  \Z[\Lambda^{*}]^{W_{J}}.
\]
\begin{proof}
Take any irreducible representation $\chi_{\mu} \in R(\hat{G})_{k}$ whose highest weight is given by $(\mu, k)$. The map
\begin{equation}
\label{iso}
\chi_{\mu} \mapsto Sk_{I}(\mu + \rho)
\end{equation}
gives the desired isomorphism. 
\end{proof}
\end{theorem}

\begin{remark}
If we identify $R(G)$ with $R(\hat{G})_{0}$, then $R(\hat{G})_{k}$ has a $R(G)$-module structure given by multiplication. On the other hand, 
 the space $\Z[\Lambda^{*}]^{W_{I}}$ is invariant under the action of $\Z[\Lambda^{*}]^{W}$ by multiplication. Therefore, the map (\ref{iso}) is also an isomorphism between $R(G) \cong \Z[\Lambda^{*}]^{W}$-modules. 
\end{remark}

Let  $\Z^{-\infty}[\Lambda^{*}]$ be the completion of $\Z[\Lambda^{*}]$, that is, the set of all integral linear functions on the lattice $\Lambda$ without requiring to have compact support. The Theorem \ref{identify theorem} can be generalized to  the following commutative diagram:
\begin{diagram}
 R^{-\infty}(\hat{G_{I}})_{k} & \rTo^{\cong}  & \Z^{-\infty}[\Lambda^{*}]^{W_{I}}  \\
 \dTo^{\mathrm{Ind}_{I}^{J}}&  & \dTo^{\mathrm{Sk}_{I}^{J}} \\
 R^{-\infty}(\hat{G_{J}})_{k} & \rTo^{\cong}  & \Z^{-\infty}[\Lambda^{*}]^{W_{J}} \\
\end{diagram}

\section{The Verlinde Ring}
Based on Meinrenken's approach \cite[Section 5]{Meinrenken09} with mild simplification, we  will compute the equivariant twisted $K$-homology in this section. 

Let $G$ be a compact, simple and simply connected Lie group. We fix a maximal torus $T$ and a positive Weyl chamber $\mathfrak{t}_{+}^{*}$. The \emph{basic inner product} on $\g$ is the unique ad-invariant inner product such that the norm of the highest root $\alpha_{\mathrm{max}} = -\alpha_{0}$ equals to $\sqrt{2}$. Using this inner product, we can identify $\g$ and $\g^{*}$. Let us fix a level $k$ and define $\mathfrak{A}_{k} \subset \Lambda^{*}$ to be  the bounded alcove cut out by the inequalities
\[
\langle \alpha_{i}^{\vee}, \cdot \rangle + k \cdot \delta_{i,0} \geq 0, \ i = 0, \dots, l.
\]
In particular, $\nu \in \mathfrak{A}_{k}$ if and only if $\frac{\nu}{k}$ lies in the fundamental alcove $\mathfrak{A}$. We observe that the weights in the shifted alcove $\mathfrak{A}_{k}^{\rho} = \mathfrak{A}_{k} + \rho$ are exactly the weights in the interior of the enlarged alcove $\mathfrak{A}_{k + h^{\vee}}$. Moreover, for any $\lambda \in \Lambda^{*}$, there exists a $w \in W_{\mathrm{aff}}$ such that 
\[
\lambda \in w \cdot \mathfrak{A}_{k + h^{\vee}} = \big\{ w \star_{k} \nu \big| \nu \in \mathfrak{A}_{k + h^{\vee}} \big\}. 
\]
The affine reflections across the faces of $\mathfrak{A}_{k + h^{\vee}}$ generate the $W_{\mathrm{aff}}$-action. We introduce a length function : $\Lambda^{*} \to \Z$.
\begin{definition} 
For any $\lambda \in \Lambda^{*}$, we define
\[
\mathrm{length}(\lambda) = \mathrm{min}\{ \mathrm{length}(w)| w \star_{k} \lambda \in \mathfrak{A}_{k + h^{\vee}}, w \in W_{\mathrm{aff}}\}.
\]
For any two weights $\lambda, \mu \in \Lambda^{*}$, if they are in the same $W_{\mathrm{aff}}$-orbit, then we define
\[
\mathrm{length}(\lambda, \mu) = \mathrm{min}\{ \mathrm{length}(w)| w \star_{k} \lambda = \mu,  w \in W_{\mathrm{aff}}\}.
\]
Otherwise, we define $\mathrm{length}(\lambda, \mu)  = \infty$. In addition, for any $\lambda$ and alcove $w_{\alpha} \cdot \mathfrak{A}_{k+h^{\vee}}$, we define their distance to be 
\[
\mathrm{min} \{\mathrm{length}(w)\big| w \star_{k} \lambda \in w_{\alpha} \cdot \mathfrak{A}_{k+h^{\vee}},  w \in W_{\mathrm{aff}} \}.
\] 
\end{definition}

For each $\lambda \in \mathfrak{A}_{k}$, we define a special element
\[
t_{\lambda} = \mathrm{exp}(\frac{\lambda+ \rho}{k + h^{\vee}}) \in T. 
\]
\begin{definition}\cite{Verlinde88}
The \emph{level k Verlinde ring} is the quotient
\[
R_{k}(G) = R(G)/I_{k}(G),
\]
where $I_{k}$ is the \emph{level k fusion ideal}:
\[
I_{k} =  \{ \chi \in R(G) \big| \chi(t_{\lambda}) = 0, \forall \ \lambda \in \mathfrak{A}_{k} \}.
\]
\end{definition}

\begin{remark}
The Verlinde ring $R_{k}(G)$ can also be defined as the fusion ring of level $k$ positive energy representations of the loop group $LG$. As a $\Z$-module, we can identify
$R_{k}(G) \cong \Z[\mathfrak{A}_{k}^{\rho} ]$. 
\end{remark}

To compute the $G$-equivariant twisted $K$-homology $KK_{\bullet}(C_{\A^{k+h^{\vee}}}(G), \C)$,  the formula (\ref{G}) realizes $\A$ geometrically as a co-simplicial Dixmier-Douady bundle, with $p$-simplices given by
\[
\coprod_{|I|= p+1}\A_{I} \to \coprod_{|I|=p+1}G/G_{I}.
\]
We obtain a spectral sequence:
\[
E^{1}_{p,q} = \bigoplus_{|I|=p+1} KK^{G}_{q}(C_{\A_{I}^{k+h^{\vee}}}(G/G_{I}), \C),
\]
with the differential $d^{1}$ given by an alternating sum:
\[
d^{1} = \sum_{r=0}^{p}(-1)^{r} \phi_{I}^{\delta_{r}I},
\]
where $\delta_{r}I = \{i_{0}, \dots, \hat{i}_{r}, \dots, i_{p} \}$ is a subset of $I$ obtained by omitting the $r$-th entry. By the Bott periodicity and the fact that 
\[
KK_{1}^{G}(C_{\A_{I}^{k+h^{\vee}}}(G/G_{I}), \C) = 0,
\] 
the $E^{1}$-sequence collapses to a single chain $(E^{1}_{\bullet}, \partial)$:
\begin{equation}
\label{chain-cplx}
0 \to E^{1}_{l} \to \dots \to E^{1}_{0} \to 0,
\end{equation}
where 
\[
E^{1}_{p} = \bigoplus_{|I|=p+1} KK^{G}_{0}(C_{\A_{I}^{k+h^{\vee}}}(G/G_{I}), \C) \cong \bigoplus_{|I|= p+1} \Z[\Lambda^{*}]^{W_{I}}
\]
In terms of weights, we can express the differential 
\[
\partial : \Z[\Lambda^{*}]^{W_{I}} \to \bigoplus_{r} \Z[\Lambda^{*}]^{W_{\delta_{r}I}}
\] 
by
\begin{equation}
\label{differential}
\partial Sk_{I}(\lambda) := \sum (-1)^{r} \cdot Sk_{\delta_{r}I}(\lambda) \in  \bigoplus_{r} \Z[\Lambda^{*}]^{W_{\delta_{r}I}}.
\end{equation}
We now pass to the $E^{2}$-sequence:
\[ 
E^{2}_{p} = \frac{\mathrm{Ker}(\partial : E^{1}_{p} \to E^{1}_{p+1})}{\mathrm{Im}(\partial : E^{1}_{p-1} \to E^{1}_{p})}.
\]
The next lemma will be the key tool in the computation.
\begin{lemma}
\label{lemma-1}
Every element in $E^{2}_{p}$ can be represented by an element in the form of  $\sum_{J} Sk_{J}(\lambda_{J})$, where $\lambda_{J}$ is scalar multiple of weights in $\mathfrak{A}_{k + h^{\vee}}$. 
\begin{proof}
Take an arbitrary element $[Sk_{I}(\lambda)] \in E^{2}_{p}$, where $\lambda$ is scalar multiple of a weight in $\mathfrak{t}_{+}^{*}$. There exists a $w \in W_{\mathrm{aff}}$ such that
\[
\mathrm{length}(w) = \mathrm{length}(\lambda), \  w \star_{k} \lambda \in  \mathfrak{A}_{k + h^{\vee}}.
\]
We can write $w =  w^{'} \cdot w_{s}$  such that
\[
 \mathrm{length}(w^{'}) = \mathrm{length}(w) -1. 
\] 
If $s \in I$, then we have that 
\[ Sk_{I}(\lambda) = - Sk_{I}(w_{s} \star_{k} \lambda) \in E^{2}_{p}, \]
and $\mathrm{length}(w_{s} \star_{k} \lambda) = \mathrm{length}(\lambda) - 1$. Otherwise, we denote 
\[
\tilde{I} = I \sqcup \{ s\} \ \mathrm{and} \ \widetilde{\delta_{r}I} =\delta_{r}I \sqcup \{s\}.
\] 
Since
\[
\partial Sk_{\tilde{I}}(\lambda) = \sum(-1)^{r} Sk_{\delta_{r}\tilde{I}}(\lambda),
\]
we have that
\[ 
[Sk_{I}(\lambda)] = \sum (-1)^{\pm} [Sk_{\widetilde{\delta_{r}I}}(\lambda))] \in E^{2}_{p}. 
\]
For every term in the right-hand side of the above equation, $\widetilde{\delta_{r}I}$ contains $s$. As in the first case, we get that 
\[
[Sk_{I}(\lambda)] =  \sum (-1)^{\pm}  \big[Sk_{\widetilde{\delta_{r}I}}(w_{s} \star_{k} \lambda) \big]  \in E^{2}_{p}. 
\]
To sum up, for any element $[Sk_{I}(\lambda)] \in E^{2}_{p}$, we can always find an equivalent element
\[ 
 \sum_{J} \big[Sk_{J}(\lambda_{J}) \big] = [Sk_{I}(\lambda) ] \in E^{2}_{p},
\]
such that $\mathrm{length}(\lambda_{J}) = \mathrm{length}(\lambda) - 1$ for all $J$. Therefore, we can prove the lemma by induction. 
\end{proof}

\begin{remark}
\label{moving}
By Lemma \ref{lemma-1}, we can move any weight $\mu \in \Lambda^{*}$ to an arbitrary alcove $w_{\alpha} \cdot \mathfrak{A}_{k + h^{\vee}}$. 
That is, for any $I \subseteq \{0, \dots, l\}$, there exists a weight $\mu^{'} \in w_{\alpha} \cdot \mathfrak{A}_{k + h^{\vee}}$ and a collection $P$ of subsets of $\{0, \dots, l\}$ such that
\[
Sk_{I}(\mu) \sim \sum_{J \in P} (-1)^{\pm} Sk_{J}(\mu^{'}) \in E^{2}_{p}. 
\]
Suppose that $z \in E^{1}_{p+1}$ realizes the equivalence, that is
\[
Sk_{I}(\mu) - \sum_{J \in P} (-1)^{\pm} Sk_{J}(\mu^{'}) = \partial z.
\]
There are two possibilities: either 
\[
\mathrm{length}(\mu, \mu^{'}) = \mathrm{length}(\mu) + \mathrm{length}(\mu^{'}),
\]
or
\[
\mathrm{length}(\mu, \mu^{'}) = \vert \mathrm{length}(\mu) - \mathrm{length}(\mu^{'}) \vert. 
\]
In the first case,  any weight whose length is less than 
\[
\mathrm{max} \big\{ \mathrm{length}(\mu), \mathrm{length}(\mu^{'}) \big\}
\] 
may appear in $z$; and in the second case $z$ involves only  those weights whose length are between $\mathrm{length}(\mu)$ and $\mathrm{length}(\mu^{'})$. This observation will play an important role in the crossed product case. 
\end{remark}

\end{lemma}

\begin{theorem}
\label{computation-2}
The homology $E^{2}_{\bullet}$ of the chain complex $E^{1}_{\bullet}$ is isomorphic to the following
\begin{eqnarray}
E^{2}_{p} \cong
\begin{cases}
0 & p = l\\
0      &  0 < p < l \\
R_{k}(G)  & p =0 \\
\end{cases}
\end{eqnarray}
\begin{proof}
\begin{enumerate}
\item
Let us start with the case when $p=l$. Take an arbitrary $[\lambda] \in E^{2}_{l}$. If
\[
\partial \lambda = \sum_{r=0}^{l} (-1)^{r}Sk_{\delta_{r}I}(\lambda) = 0,
\]
then each term $Sk_{\delta_{r}I}(\lambda) = 0$. This implies
\[
w_{r} \star_{k}\lambda = \lambda, \ r = 0, \dots, l.
\]
Hence, $\lambda$ is invariant under the affine Weyl group action.  But $\lambda$ has compact support, it must be 0.

\item
In the case when $p =0$, we have that
\[
E^{2}_{0} = \bigoplus_{i=0}^{l} \Z[\Lambda^{*}]^{W_{i}}.
\]
Take an element $Sk_{i}(\lambda) \in E^{2}_{0}$, where  $\lambda$ is scalar multiple of weights in $\mathfrak{A}_{k+h^{\vee}}$. For any $i, j$, the two elements $Sk_{i}(\lambda), Sk_{j}(\lambda)$ are differ by the boundary of $Sk_{ij}(\lambda) \in E^{2}_{1}$. Hence, 
\[
[Sk_{i}(y)] = [Sk_{j}(y)]  \in E^{2}_{0}.
\]  
On the other hand, if $\lambda \in \mathfrak{A}_{k+h^{\vee}} \setminus \mathfrak{A}_{k}^{\rho} $, then there exists an element $w_{j}$ in the Weyl group such that $w_{j} \star_{k} \lambda = 0$. Hence,
\[
[Sk_{i}(\lambda)] = [Sk_{j}(\lambda)] = [-Sk_{j}(w_{j} \star_{k} \lambda) ] = 0 \in E^{2}_{0}.
\]
By Lemma \ref{lemma-1},  $E^{2}_{0}$ is identical with $\Z[\mathfrak{A}_{k}^{\rho}] \cong R_{k}(G)$.

\item
When $ 0 < p <l$, using Lemma \ref{lemma-1} again, let us take $ x \in E^{2}_{p}$:
\[ 
 x =  \sum_{I \in \mathscr{A}}  Sk_{I}(\lambda) \in E^{2}_{p},
\]
where $\lambda$ is scalar multiple of a weight in $\mathfrak{A}_{k+h^{\vee}}$. We choose a $s \in \{ 0, \dots, l\}$ such that the two sets
\[
\mathscr{B} = \{ I \in \mathscr{A} \big| s \in I\}  \ \mathrm{and} \ \mathscr{C} = \{ I \in \mathscr{A}  \big| s \notin I\}
\]
are both non-empty. For simplicity, we denote 
\[
\mathscr{D} = \{ I \subset \{0, \dots, p\} \big| I\sqcup \{s\} \in \mathscr{B}\}.
\]
By formula (\ref{differential}), we have that
\begin{equation}
\begin{aligned}
\partial x & = \sum_{I \in \mathscr{D}}\sum_{r} (-1)^{\pm} Sk_{\delta_{r}I\sqcup \{s\}}(\lambda) \\
& + \sum_{I \in \mathscr{D}} (-1)^{\pm} Sk_{I}(\lambda) + \sum_{I \in \mathscr{C}}\sum_{r} (-1)^{\pm} Sk_{\delta_{r}I}(\lambda) \in E^{2}_{p}
\end{aligned}
\end{equation}
If $\partial x = 0$, then the last two terms in above equation must be canceled since they do not contain $s$. This implies that 
\[
\sum_{I \in \mathscr{D}} (-1)^{\pm} Sk_{I}(\lambda)
\]
is the boundary of $\sum_{I \in \mathscr{C}} (-1)^{\pm} Sk_{I}(\lambda) \in E^{1}_{p}$. As a result, one can deduce that 
\[
x = \sum_{I \in \mathscr{A}} Sk_{I}(\lambda) \in E^{2}_{p}
\] 
is the boundary of 
\[
\sum_{I \in \mathscr{C}}  Sk_{I \sqcup \{s\}}(\lambda) \in E^{1}_{p+1}. 
\]
 \end{enumerate}
\end{proof}
\end{theorem}

By Theorem \ref{computation-2}, we can recover Freed-Hopkins-Teleman's result for a special case:
\begin{corollary}
For a compact, simple and simply connected Lie group $G$, we have that
\[
KK^{G}_{\bullet}(C_{\A}(G), \C) \cong R_{k}(G). 
\]
\end{corollary}

\begin{remark}
By the Poincare duality for the twisted $K$-homology and twisted $K$-theory obtained by J.-L. Tu \cite{Tu09}, the twisted $K$-theory
\[
K_{G}(G, \A^{k+h^{\vee}}) = KK^{G}(\C, C_{\A^{k+h^{\vee}}}(G))
\]
is isomorphic to $R_{k}(G)$.
\end{remark}


\section{The Formal Verlinde Module}
 We now extend the techniques in last section to compute the twisted $K$-homology of crossed product:
\[
KK_{\bullet}(C^{*}(G, C_{\mathcal{A}^{k+h^{\vee}}}(G)), \C).
\]
This section is the heart of this paper. 

As before, we construct a spectral sequence 
\[
E^{1}_{p,q} = \bigoplus_{|I|=p+1} KK_{q}(C^{*}(G, C_{\A_{I}^{k+h^{\vee}}}(G/G_{I})), \C),
\]
which collapses to a single chain $(E^{1}_{\bullet}, \partial)$.  In terms of weights, the differential $\partial$ is given by
\[
\partial Sk_{I}(\lambda) := \sum (-1)^{r} \cdot Sk_{\delta_{r}I}(\lambda) \in \bigoplus \Z^{-\infty}[\Lambda^{*}]^{W_{\delta_{r}I}}.
\]
Note that $\lambda \in \Z^{-\infty}[\Lambda^{*}]$ may not have compact support. The computation for the crossed product case is similar to the classical case but slightly more elaborate. In the classical case, we can always simplify the computation by moving all the weights to the alcove $\mathfrak{A}_{k+h^{\vee}}$ using Lemma \ref{lemma-1}. However, it is no longer true for the crossed product case.  

For any element $[x] \in E^{2}_{p}$, suppose that it is represented by
\begin{equation}
\label{x}
x = \sum_{\beta \in P} Sk_{I_{\beta}}(\lambda_{\beta}) \in E^{2}_{p}, 
\end{equation}
where $\lambda_{\beta}$ are scalar multiple of weights, and $P$ is a collection of subsets of $\{0, \dots, l\}$. Take a sequence $\{w_{\alpha_{i}} \} \in W_{\mathrm{aff}}$. We obtain a collection of alcoves $\big\{w_{\alpha_{i}} \cdot \mathfrak{A}_{k+h^{\vee}} \big\}$. 

\begin{definition}
We say that the collection of alcoves $\big\{w_{\alpha_{i}} \cdot \mathfrak{A}_{k+h^{\vee}}^{\rho} \big\}$ and $x$ are \emph{separated by a finite distance} if there exists a constant $K$ such that for all $\lambda_{\beta}$ in (\ref{x}), there is an alcove $w_{\alpha_{i}} \cdot \mathfrak{A}_{k+h^{\vee}}$ whose distance to $\lambda_{\beta}$ is less than $K$. 
\end{definition}

For example, $Sk_{I}(\lambda)$ and  $\mathfrak{A}_{k+h^{\vee}}$  are separated by a finite distance if and only if $\lambda$ has compact support. 

\begin{lemma}
\label{lemma-3}
If a collection of alcoves $\big\{w_{\alpha_{i}} \cdot \mathfrak{A}_{k+h^{\vee}} \big\}$ and $x$ are separated by a finite distance, then we can find an element
\[
 y = \sum_{i}\sum_{I \in P_{i}}Sk_{I}(\mu_{i,I}) \in E^{2}_{p}
\]
such that $[x] = [y] \in E^{2}_{p}$, and $\mu_{i,I}$ are scalar multiple of weights in $w_{\alpha_{i}} \cdot \mathfrak{A}_{k+h^{\vee}}^{\rho}$  for any $I \in P_{i}$. 
\begin{proof}
If $x$ does not have compact support, we are not able to move all the weights $\lambda_{\beta}$ in (\ref{x}) to the alcove $\mathfrak{A}_{k+h^{\vee}}$ simultaneously. Otherwise, suppose that $z$ realizes this operation. By Remark \ref{moving}, the element $z$ is not well-defined since it may contain coefficients equaling to infinity. We will fix it in the following way. 

We re-arrange the expression of $x$ in  (\ref{x}):
\begin{equation}
\label{x-1}
x  = \sum_{i=1}^{\infty} \sum_{I \in Q_{i}} Sk_{I}(\eta_{i,I}), 
\end{equation}
where $\eta_{i,I}$ are scalar multiple of  weights, whose distances to the alcove $w_{\alpha_{i}}\cdot \mathfrak{A}_{k+h^{\vee}}$ are less than $K \gg 0$. We can move 
all the $\eta_{i,I}$ to the alcove $w_{\alpha_{i}}\cdot \mathfrak{A}_{k+h^{\vee}}$. Suppose that $z_{i} \in E^{1}_{p+1}$ realizes the moving operation. That is,
\[
 \partial z_{i} = \sum_{I \in Q_{i}} Sk_{I}(\eta_{i,I}) - \sum_{J \in P_{i}}Sk_{J}(\mu_{i,J})\in E^{1}_{p},
\] 
where $ \mu_{i,J}$ are scalar multiple of  weights in $w_{\alpha_{i}}\cdot \mathfrak{A}_{k+h^{\vee}}$. By Remark \ref{moving} and the assumption that the distances between weights and alcoves are bounded by $K$, the infinite sum 
\[
z = \sum_{i}z_{i} \in E^{1}_{p+1}, 
\]
is well-defined. If we define 
\[
y = \sum_{i}  \sum_{J \in P_{i}} Sk_{J}(\mu_{i,J}) \in E^{1}_{p},
\] 
then $x$ and $y$ is differed by the boundary of $z$. 
\end{proof}

\end{lemma}

\begin{proposition}
\label{p1}
At the $l = \mathrm{dim}(G)$ level, we have 
\[
E^{2}_{l} \cong \Z^{-\infty}[\Lambda^{*}]^{W_{\mathrm{aff}}},
\]
where $\Z^{-\infty}[\Lambda^{*}]^{W_{\mathrm{aff}}}$ is the invariant part of $\Z^{-\infty}[\Lambda^{*}]$ for the $W_{\mathrm{aff}}$-action.
\begin{proof}
Take an arbitrary $ [\lambda] \in E^{2}_{l}$. The condition
\[
\partial \lambda = \sum_{r=0}^{l} (-1)^{r}Sk_{\delta_{r}I}(\lambda) = 0
\]
implies $\lambda$ is invariant under the $W_{\mathrm{aff}}$-action. 
\end{proof}
\end{proposition}

\begin{proposition}
\label{p2}
For any $0 < p <l$, we have 
\[
E^{2}_{p} = 0.
\]
\begin{proof}
We fix a  constant $C \gg 0$ and choose a sequence $\{w_{\alpha_{i}} \}_{i=1}^{\infty}$ in $W_{\mathrm{aff}}$ such that for any $i \neq j$
\[
 \vert \mathrm{length}(w_{\alpha_{i}}) -   \mathrm{length}(w_{\alpha_{j}}) \vert \geq C. 
\]
In this case, for any two weights $\lambda_{\alpha_{i}} \in w_{\alpha_{i}}\cdot \mathfrak{A}_{k+h^{\vee}}$ and $\lambda_{\alpha_{j}} \in w_{\alpha_{j}}\cdot \mathfrak{A}_{k+h^{\vee}}$, 
\[
\partial Sk_{I}(\lambda_{\alpha_{i}}), \partial Sk_{J}(\lambda_{\alpha_{j}}) \in E^{2}_{p}, 
\]
do not share common terms for arbitrary $I, J$. 

For any element $[x] \in E^{2}_{p}$, using Lemma \ref{lemma-3}, we can assume that
\[
x = \sum_{i=1}^{\infty} \sum_{I \in P_{i}} Sk_{I}(\lambda_{i,I}), 
\]
where $\lambda_{i, I}$ are scalar multiple of weights in $w_{\alpha_{i}} \cdot \mathfrak{A}_{k+h^{\vee}}$, for  all $I \in P_{i}$. By the assumption that the  alcoves are separated, the fact $\partial x = 0$ implies that $\partial \sum_{I \in P_{i}} Sk_{I}(\lambda_{i, I}) = 0$ individually for each $i$. From here,  the argument in Theorem \ref{computation-2} may be applied verbatim. 
\end{proof}
\end{proposition}

At last, when $p=0$, we have that 
\[
E^{1}_{0} = \bigoplus_{i=0}^{l} \Z^{-\infty}[\Lambda^{*}]^{W_{i}}.
\] 
For any element $x \in E^{1}_{0}$, we  express it by
\[
x = \sum_{i=0}^{l}Sk_{i}(\lambda_{i}),
\]
where $\lambda_{i} \in \in \Z^{-\infty}[\Lambda^{*}]$. In $E^{2}_{0}$, we have that 
\[
[Sk_{0}(\lambda)] = [Sk_{i}(\lambda)]  = [-Sk_{0}(w_{j} \star_{k}\lambda)] \in E^{2}_{0},
\]
for any $1\leq i,j \leq l$. Note that for any element $\lambda \in \Z^{-\infty}[\Lambda^{*}]^{W_{\mathrm{aff}}}$, 
\[
Sk_{i}(\lambda) = 0, \ \mathrm{for} \ i = 0, \dots, l.
\]
Hence, equivalence classes in $E^{2}_{0}$ are generated by $[Sk_{0}(\lambda)]$, where 
\[
\lambda \in \breve{\Z}^{-\infty}[\Lambda^{*}] = \Z^{-\infty}[\Lambda^{*}]^{W} \setminus \Z^{-\infty}[\Lambda^{*}]^{W_{\mathrm{aff}}}.
\] 
Here, the set  $\Z^{-\infty}[\Lambda^{*}]^{W}$ is the anti-invariant part of  $\Z^{-\infty}[\Lambda^{*}]$ for the Weyl group action. That is, $v \in \Z^{-\infty}[\Lambda^{*}]^{W}$ if and only if for all $w \in W$
\[
w \star v = w \cdot (v + \rho) - \rho = (-1)^{\mathrm{length}(w)}\cdot v.
\] 
The set $\Z^{-\infty}[\Lambda^{*}]^{W}$ has a $\Z[\Lambda^{*}]^{W}$-module structure while $R^{-\infty}(G)$ has a $R(G)$-module structure as well. Thus, by identifying $R(G)$ and $\Z[\Lambda^{*}]^{W}$, $\Z^{-\infty}[\Lambda^{*}]^{W}$ is isomorphic to $R^{-\infty}(G)$ as $R(G)$-modules. 
For simplicity, we will not distinguish them for the rest of this section. 
\begin{lemma}
\label{lemma-2}
For any element $\lambda \in \breve{\Z}^{-\infty}[\Lambda^{*}]$, if 
\[
[Sk_{0}(\lambda)] = 0 \in E^{2}_{0},
\] 
then, as an element in $R^{-\infty}(G)$
\[
\lambda \in I_{k}(G) \cdot R^{-\infty}(G). 
\]
\begin{proof}
First of all, if $\lambda$ has compact support, by the conclusion in the classical case, we have that $\lambda \in I_{k}(G)$. 

In general, we fix two constants $K  \gg C \gg 0$. By Lemma \ref{lemma-3}, there exists an element $\mu \in \breve{\Z}^{-\infty}[\Lambda^{*}]$ such that
\[
[Sk_{0}(\lambda)] = [Sk_{0}(\mu)] \in E^{2}_{0},
\]
and the support of $\mu$ is discrete. To be more precise, we can rewrite
\[
\mu = \sum_{i=0}^{\infty} \sum_{\alpha \in P_{i}} \mu_{\alpha}^{i} \in \breve{\Z}^{-\infty}[\Lambda^{*}],
\] 
where $\mu_{\alpha}^{i}$ are scalar multiple of weights and $\vert \mathrm{length}(\mu_{\alpha}^{i}) - i \cdot K\vert \leq C$. Because $[Sk_{0}(\mu)] = 0\in E^{2}_{0}$ and the assumptions on $\mu$, the element 
\[
\sum_{\alpha \in P_{i}} Sk_{0}( \mu_{\alpha}^{i}),
\] 
which has compact support, must be  boundary as well. From last section, we know that $\sum_{\alpha \in P_{i}} \mu_{\alpha}^{i} \in I_{k}(G)$. By Remark \ref{moving}, the sum 
\[
\mu = \sum_{i=0}^{\infty}\mu_{i} \in I_{k}(G)\cdot R^{-\infty}(G).
\]

\end{proof}
\end{lemma}

\begin{proposition}
\label{p3}
At $p =0$ level, we have that
\[
E^{2}_{0} \cong \breve{\Z}^{-\infty}[\Lambda^{*}] \otimes_{R(G)} R_{k}(G). 
\]
\begin{proof}
This proposition follows from Theorem \ref{computation-2} and Lemma \ref{lemma-2}. 
\end{proof}
\end{proposition}

By Proposition \ref{p1}-\ref{p3}, we can now state the main theorem of this paper. 

\begin{theorem}
For any compact, simple, simply connected Lie group $G$, and $\A \to G$ the $G$-equivariant Dixmier-Douady bundle corresponding to $[1] \in H_{G}^{3}(G) \cong \Z$, we have the identity
\[
KK_{\bullet}(C^{*}(G,  C_{\A^{k + h^{\vee}}}(G)), \C) \cong R^{-\infty}(G) \otimes_{R(G)} R_{k}(G).
\]
\end{theorem}

\section{Quantization and K-Homology}
In this section we discuss some application of the formal Verlinde module in the quantization of quasi-Hamiltonian spaces \cite{Alekseev98}, which is our original motivation. 

Let $G$ be a compact, connected Lie group, and $(M, \omega)$ be a compact symplectic $G$-manifold. The symplectic 2-form $\omega$ determines a $G$-equivariant Spin$^{c}$-structure, thus a $G$-equivariant Morita morphism:
\[
(\phi, S_{M}) : (M, \mathrm{Cliff}(TM)) \to (\pt, \C).
\]
The morphism induces  a push-forward map in  $K$-homology:
\begin{equation}
\label{push-forward}
\phi_{S_{M}} :KK_{0}^{G}(C_{\tau}(M), \C) \to KK_{0}^{G}(\C,\C) \cong R(G),
\end{equation}
where $C_{\tau}(M)$ is the $C^{*}$-algebra of continuous sections of the complex Clifford bundle $\mathrm{Cliff}(TM)$. In addition, the map (\ref{push-forward}) can be extended from $G$-equivariant $K$-homology to $K$-homology of crossed products of $C^{*}$-algebras:
\[
\hat{\phi}_{S_{M}} :KK_{0}(C^{*}(G, C_{\tau}(M)), \C) \to KK_{0}(C^{*}(G),\C) \cong R^{-\infty}(G). 
\]

On the other hand, Kasparov \cite{Kasparov88} defined a distinguished element 
\[
[d_{M}] \in KK_{0}^{G}(C_{\tau}(M), \C),
\] 
the \emph{Dirac element}, using the de Rham differential operator $d$ and its dual $d^{*}$. With this Dirac element, Meinrenken \cite{M12} defined the quantization of the Hamiltonian $G$-space $(M, \omega)$ as the push-forward
\[
Q(M, \omega) := \phi_{S_{M}}([d_{M}]) \in KK_{0}^{G}(\C,\C) \cong R(G).
\]
This reformulation of  quantization no longer mentions Spin$^{c}$-Dirac operators which are blur in the quasi-Hamiltonian case. 

Following an idea of Witten \cite{Witten82}, we can perturb the de Rham differential operator $d$ by some $G$-invariant smooth function $f$ as follows:
\[
d_{f} = e^{-f} \cdot d \cdot e^{f}, \ d_{f}^{*} = e^{f} \cdot d^{*} \cdot e^{-f}. 
\]
Let us define a deformed Dirac element $[d_{M, f}]$. If the set of critical points of $f$ can be separated by a collection of disjoint open sets $U_{\alpha}$, then we can break 
\[
[d_{M, f}] = \sum_{\alpha} [d_{U_{\alpha},f}]
\]
in the $K$-homology group $KK_{0}(C^{*}(G, C_{\tau}(M)), \C)$. The quantization of each piece
\[
\hat{\phi}_{S_{M}}([d_{U_{\alpha, f}}]) \in KK_{0}(C^{*}(G), \C) \cong R^{-\infty}(G),
\] 
can be interpreted as the index of some transversally elliptic operators \cite{Atiyah74, Kasparov14} associated the open symplectic manifolds $U_{\alpha}$ and function $f$. We obtain a localization formula:
\begin{equation}
\label{local}
Q(M, \omega) = \hat{\phi}_{S_{M}}([d_{M}]) = \sum_{\alpha} \hat{\phi}_{S_{M}}([d_{U_{\alpha, f}}])\in KK_{0}(C^{*}(G), \C) \cong R^{-\infty}(G). 
\end{equation}

Let us now consider  a quasi-Hamiltonian space $M$ with a group-valued moment map $\phi$ \cite{Alekseev98}. In general, there is no Spin$^{c}$-structure on $M$. However, we can define a twisted Spin$^{c}$-structure given by a Morita morphism:
\[
(\phi, \mathcal{E}) : (M, \mathrm{Cliff}(TM)) \to (G, \A^{k + h^{\vee}}),
\]
which induces a map 
\[
\phi_{*} : KK^{G}_{0}(C_{\tau}(M), \C) \to KK^{G}_{0}(C_{\A^{k+h^{\vee}}}(G), \C).
\]
Meinrenken \cite{M12} defined the quantization of the quasi-Hamiltonian space $(M, \phi)$ as the push-forward of the Dirac element:
\[
Q(M) = \phi_{*}([d_{M}]) \in KK^{G}_{0}(C_{\A^{k+h^{\vee}}}(G), \C) \cong R_{k}(G). 
\]

It is natural to ask if we can generalize the localization formula (\ref{local}) to the quasi-Hamiltonian case. As before, we choose a $G$-invariant function $f$ and break the Dirac element
\[
[d_{M, f}] = \sum_{\alpha} [d_{U_{\alpha},f}] \in KK_{0}(C^{*}(G, C_{\tau}(M)),\C).
\]
By the push-forward map, we have that
\[
\hat{\phi}_{*}([d_{U_{\alpha},f}]) \in KK_{0}(C^{*}(G, C_{\A^{k+h^{\vee}}}(G)) , \C) \cong R^{-\infty}(G) \otimes_{R(G)} R_{k}(G). 
\] 
This suggests that the formal Verlinde module
\[
R^{-\infty}(G) \otimes_{R(G)} R_{k}(G)
\]
could be interpreted as quantization of some non-compact quasi-Hamiltonian spaces (of course we have to impose some suitable properness  assumptions to make sense of it).

\bibliographystyle{alpha}
\bibliography{mybib}

\end{document}